\documentclass[a4paper,fleqn, sort&compress]{cas-sc}

\usepackage[numbers]{natbib}
\usepackage[nameinlink,capitalise,noabbrev]{cleveref}

\usepackage{amsmath, amssymb, amsthm, bm}
\usepackage{mathtools}
\newtheorem{theorem}{Theorem}[section]
\newtheorem{lemma}[theorem]{Lemma}
\newtheorem{proposition}[theorem]{Proposition}

\theoremstyle{remark}

\theoremstyle{definition}

\usepackage{changepage}

\usepackage[utf8]{inputenc}
\usepackage[T1]{fontenc}
\usepackage{longtable}
\usepackage{xfrac}
\usepackage{graphicx}
\usepackage{amsmath,amssymb,amsfonts}
\usepackage{csquotes}
\usepackage{subcaption}
\usepackage{caption}
\usepackage{capt-of}
\usepackage{algorithm}
\usepackage{algpseudocode}
\usepackage{tabularx}
\usepackage{amsmath}
\usepackage{amssymb}
\usepackage{mathtools}

\usepackage{booktabs} 
\usepackage{tabularray}
\UseTblrLibrary{booktabs}

\usepackage{placeins}

\usepackage{enumitem}
\usepackage{amsthm}%

\PassOptionsToPackage{hyphens}{url} 
\usepackage{hyperref}
\hypersetup{breaklinks=true}  
\def\tsc#1{\csdef{#1}{\textsc{\lowercase{#1}}\xspace}}
\tsc{WGM}
\tsc{QE}
\tsc{EP}
\tsc{PMS}
\tsc{BEC}
\tsc{DE}

\begin{document}
\let\WriteBookmarks\relax
\def\floatpagepagefraction{1}
\def\textpagefraction{.001}

\shorttitle{One-Step Information Boundaries in Finite-Volume Coarse Graining}

\shortauthors{A. Polemitis, N. Christakis and D. Drikakis}

\title [mode = title]{Initial-Trace Insufficiency and Entropy-Accounting Separation in Finite-Volume Coarse Graining}                      



%
\author[1]{Antonis Polemitis}[orcid=0000-0003-3188-8746]
\cormark[1]
\ead{polemitis.an@unic.ac.cy}

\author[2]{Nicholas Christakis}[
                        orcid=0000-0001-9212-2336]

\cormark[2]


\ead{christakis.n@unic.ac.cy}



\affiliation[1]{UNIC Evolve, University of Nicosia, Nicosia, CY-2417, Cyprus}

\affiliation[2]{Institute for Advanced Modelling and Simulation, University of Nicosia, Nicosia, CY-2417, Cyprus}

\author[2]{Dimitris Drikakis}[orcid=0000-0002-3300-7669]
\cormark[1]
\ead{drikakis.d@unic.ac.cy}



\cortext[cor1]{Corresponding author}
\cortext[cor2]{Principal corresponding author}



\begin{abstract}
Coarse finite-volume information does not automatically define an exact
pathwise closure for the restricted evolution of a fine-grid realization.  We
study two one-step information boundaries that are distinct from recursive
predictive memory.  First, we show that parent averages together with
complete initial parent-interface trace pairs are sufficient for one
forward-Euler finite-volume transition with a deterministic two-state
numerical flux, but they do not generally determine a completed
second-order strong-stability-preserving Runge--Kutta, equivalently Heun,
transition.  An explicit scalar upwind collision shows that hidden interior
data can affect later stage traces and hence the completed parent update,
even when initial traces and parent averages agree.

Second, we show that exact conservative one-step closure is not the same as
entropy-accounting closure.  A parent entropy balance obtained from a
cellwise numerical entropy balance contains a subcell Jensen gap, a realized
entropy-flux difference, and residual dissipation.  For a declared
Burgers--Rusanov conservative flux and a declared Rusanov-type numerical
entropy flux, we construct a sign-indexed finite family in which all branches
have identical parent averages, identical conservative interface registers,
and identical next parent averages, while the entropy-interface register
distinguishes every branch.  Under the uniform branch law, the conservative
channel leaves nonzero deterministic squared-error floors for both the raw
entropy register and its cyclic divergence.  The results give a concise
finite-volume information-sufficiency test for multistage interface data and
entropy-side accounting.
\end{abstract}


\begin{highlights}
\item Initial traces close forward Euler but can fail for Heun updates.
\item Realized stage-integrated fluxes close completed finite-volume transitions.
\item Conservative registers do not determine entropy-accounting information.
\item A Burgers--Rusanov sign family gives exact entropy-blindness error floors.
\end{highlights}


\maketitle

\section{Introduction}
\label{sec:introduction}

Finite-volume methods evolve cell averages through numerical fluxes at cell
interfaces \cite{LeVeque2002}.  In coarse/fine coupling, adaptive mesh
refinement, and reduced modelling, one often asks whether a declared coarse
information channel is sufficient to reproduce the restricted result of a
fine-grid computation \cite{Berger1987,BergerColella1989,
McCorquodaleColella2011,OsherSanders1983,ConstantinescuSandu2007}.  Parent
averages alone generally do not determine subcell arrangement, and subcell
arrangement can affect later interface fluxes.

A companion paper studies the recursive version of this issue for a periodic
scalar upwind hierarchy: it determines exact finite-horizon memory,
conditioning, and dissipative decay for future parent-average prediction
\cite{PolemitisChristakisDrikakis2026FVerasure}.  The present paper asks a
different one-step question: which declared interface and entropy-side data
close a completed finite-volume transition, and which proposed channels do
not?

The one-step conservative baseline is classical.  If a completed fine
finite-volume step has shared conservative face fluxes, then summing child
cell balances over a parent cancels internal fluxes and leaves only the
parent-interface flux difference.  This is the algebra behind refluxing and
conservative mesh-interface communication \cite{Berger1987,
BergerColella1989}.  The realized integrated interface flux is therefore a
canonical ex post sufficient message for reconstructing the restricted result
of an already completed fine step.  It is not, by itself, a cheap autonomous
predictor before the fine work is done.

The new points here concern two boundaries of this one-step baseline.  The
first is a multistage boundary.  For forward Euler, complete ordered initial
traces at a parent interface determine the parent-interface numerical flux.
For a multistage method, later stage traces can depend on hidden interior
fine data. We give an exact same-data collision for the second-order
strong-stability-preserving Runge--Kutta method, denoted SSPRK(2,2), also
known as the explicit trapezoidal or Heun method in this linear autonomous setting.

The second boundary is an entropy-accounting boundary.  Entropy-stable and
entropy-accounting finite-volume methods require information beyond the
conservative update \cite{Tadmor1987,Tadmor2003,FjordholmMishraTadmor2012,
Lozano2018,RanochaSayyariDalcinParsaniKetcheson2020}.  A parent entropy
balance involves a subcell Jensen gap, a realized entropy-flux difference,
and residual dissipation.  We construct an exact Burgers--Rusanov family
where the conservative channel is identical across branches while the
entropy-side register is not.

These questions are also relevant to data-driven discretizations and closure
models \cite{BarSinai:2019datadriven,KochkovEtAl2021,
AgdesteinSanderse2025,SanderseEtAl2025}.  Exact unresolved finite-volume
terms can be written down for declared filters and discretizations
\cite{AgdesteinVerstappenSanderse2026}; the present paper instead asks
whether a proposed retained channel is sufficient for a target transition or
entropy account to be single valued.

\paragraph{Contributions.}
The factorization criterion and flux telescoping identity used below are
classical.  The contributions of this paper are:
\begin{enumerate}
\item a one-step information-channel formulation for finite-volume
coarse-graining contracts;
\item a proof that initial parent-interface traces close forward Euler but
can fail for SSPRK(2,2)/Heun, even when parent averages also agree;
\item an exact Burgers--Rusanov family showing that conservative one-step
closure does not determine the augmented parent entropy account.
\end{enumerate}

\paragraph{Organization.}
\Cref{sec:baseline} states the factorization and flux-register baselines in
minimal form.  \Cref{sec:heun} proves the initial-trace sufficiency and
SSPRK(2,2)/Heun insufficiency results.  \Cref{sec:entropy} derives the
parent entropy account and the entropy-blind finite family.
\Cref{sec:discussion} discusses implications and limitations, and
\Cref{sec:conclusions} concludes.

\section{One-Step Closure Baseline}
\label{sec:baseline}

Let \(\mathcal S\) be a set of admissible fine states, let
\[
  T:\mathcal S\to\mathcal C
\]
be a declared coarse information map, and let
\[
  \mathcal H:\mathcal S\to\mathcal Y
\]
be the target restricted output.

\begin{lemma}[Fiber factorization]
\label{lem:factorization}
There exists a unique map \(G:T(\mathcal S)\to\mathcal Y\) such that
\[
  \mathcal H=G\circ T
\]
if and only if
\[
  T(x)=T(x')\quad\Longrightarrow\quad \mathcal H(x)=\mathcal H(x')
  \qquad (x,x'\in\mathcal S).
\]
\end{lemma}

\begin{proof}
If \(\mathcal H=G\circ T\), equal \(T\)-values imply equal
\(\mathcal H\)-values.  Conversely, if \(\mathcal H\) is constant on fibers
of \(T\), define \(G(t)\) as the common value of \(\mathcal H\) on
\(T^{-1}(t)\).  This is well defined on \(T(\mathcal S)\) and is unique
there.
\end{proof}

In a one-dimensional finite-volume hierarchy, a completed conservative fine
step has the form
\[
  U_i^{n+1}
  =
  U_i^n
  -
  \frac{\Phi_{i+1/2}-\Phi_{i-1/2}}{h_i},
\]
where \(\Phi_{i+1/2}\) is the realized time-integrated fine-face flux.  If
parent \(K\) contains the fine cells \(i\subset K\), and
\[
  H_K=\sum_{i\subset K}h_i,
  \qquad
  (RU)_K=\frac1{H_K}\sum_{i\subset K}h_iU_i,
\]
then summing the child balances gives
\begin{equation}
  (RU^{n+1})_K
  =
  (RU^n)_K
  -
  \frac{Q_K-Q_{K-1}}{H_K},
  \label{eq:flux-register-update}
\end{equation}
where \(Q_K\) is the realized integrated flux through the parent interface
\(K+1/2\).  On a periodic parent grid, only the cyclic differences
\(Q_K-Q_{K-1}\) affect the parent update; adding one spatially constant flux
register changes no parent average.  This elementary gauge is the one-step
counterpart of the periodic gauge used in recursive predictive-memory
analysis \cite{PolemitisChristakisDrikakis2026FVerasure}.

Equation \eqref{eq:flux-register-update} is used here only as the classical
baseline: the realized integrated flux register is sufficient for the result
of a completed step.  The next section asks when a simpler initial-trace
channel can replace that register.

\section{Initial Traces and Multistage Flux Registers}
\label{sec:heun}

Throughout this section the mesh, time step, flux parameters, and parent
partition are fixed, there are no source or nonconservative terms, and all
boundary fluxes are either periodic or externally supplied.  A trace at a
parent interface means the actual ordered one-sided state pair passed to the
deterministic numerical flux.

\begin{proposition}[Forward-Euler trace sufficiency]
\label{prop:forward-euler-trace}
Consider a one-dimensional finite-volume forward-Euler step
\[
  U_i^{n+1}
  =
  U_i^n
  -
  \frac{\Delta t}{h_i}
  \left[
  \widehat F(U_i^n,U_{i+1}^n)
  -
  \widehat F(U_{i-1}^n,U_i^n)
  \right],
\]
where \(\widehat F\) is a deterministic two-state numerical flux with known
parameters.  If the declared information consists of the current parent
averages and the complete ordered initial trace pair at every parent
interface, then the next parent averages are determined exactly.
\end{proposition}

\begin{proof}
Summing the fine updates over a parent cancels all internal fine-face fluxes.
The only remaining terms are the two parent-interface numerical fluxes
\[
  \widehat F(U_{K,r-1}^n,U_{K+1,0}^n)
  \quad\text{and}\quad
  \widehat F(U_{K-1,r-1}^n,U_{K,0}^n).
\]
These are determined by the ordered initial traces.  The corresponding
integrated forward-Euler registers are \(\Delta t\) times these numerical
fluxes.  Together with the current parent averages, they determine the next
parent averages through \eqref{eq:flux-register-update}.
\end{proof}

This forward-Euler sufficiency is special to a one-stage update.  In a
multistage method, later stage traces may depend on subcell values that are
not visible in the initial parent-interface traces.

Let \(A_\lambda\) be the periodic scalar positive-upwind forward-Euler
matrix,
\[
  (A_\lambda u)_i=(1-\lambda)u_i+\lambda u_{i-1},
  \qquad 0<\lambda\le1,
\]
with periodic indexing.  For this linear autonomous update, the
SSPRK(2,2), equivalently explicit trapezoidal or Heun, completed-step
operator is
\[
  M_\lambda
  =
  I+(A_\lambda-I)+\frac12(A_\lambda-I)^2
  =
  \frac12(I+A_\lambda^2).
\]

\begin{theorem}[Initial traces do not close SSPRK(2,2)/Heun]
\label{thm:heuntrace}
Consider \(P=2\) parents, \(r=4\) fine cells per parent, and periodic scalar
positive-speed upwinding.  Let
\[
  x^+=(0,-1,1,0\mid 0,0,0,0),
  \qquad
  x^-=(0,1,-1,0\mid 0,0,0,0),
\]
where the vertical bar separates the two parents.  The two states have
identical parent averages and identical initial parent-interface traces.
However,
\[
  R M_\lambda x^+
  =
  \left(-\frac{\lambda^2}{8},\frac{\lambda^2}{8}\right),
  \qquad
  R M_\lambda x^-
  =
  \left(\frac{\lambda^2}{8},-\frac{\lambda^2}{8}\right).
\]
Thus parent averages and initial parent-interface traces do not determine
the completed SSPRK(2,2)/Heun parent update.
\end{theorem}

\begin{proof}
Both states have zero parent averages, and the ordered trace pairs at both
parent interfaces are \((0,0)\) in both states.  Let
\[
  v=x^+-x^-=(0,-2,2,0\mid 0,0,0,0).
\]
Then \(Rv=0\).  If \(S\) denotes the periodic shift \((Su)_i=u_{i-1}\), then
\[
  A_\lambda^2
  =
  (1-\lambda)^2I
  +
  2\lambda(1-\lambda)S
  +
  \lambda^2S^2.
\]
Only the \(S^2\) term moves part of \(v\) across the parent interface.
Direct summation over the two parents gives
\[
  R M_\lambda v
  =
  \left(-\frac{\lambda^2}{4},\frac{\lambda^2}{4}\right).
\]
Since \(x^+=v/2\) and \(x^-=-v/2\), the displayed values follow.
\end{proof}

The obstruction is not that stage-integrated fluxes fail.  In a standard
Butcher-form finite-volume Runge--Kutta step, the corresponding sufficient
ex post register at parent interface \(K+1/2\) is
\[
  Q_K^{\mathrm{RK}}
  =
  \Delta t
  \sum_s b_s\,
  \widehat F\left(
    U_{K,r-1}^{(s)},U_{K+1,0}^{(s)}
  \right),
\]
with the declared stage values \(U^{(s)}\) and weights \(b_s\).  These
realized stage-integrated interface fluxes still reconstruct the restricted
result of the completed fine step.  The point is narrower and sharper: the
initial trace channel that suffices for forward Euler does not factor the
completed multistage transition.

\section{Entropy-Accounting Separation}
\label{sec:entropy}

Conservative one-step closure determines next parent averages.  It does not
determine the augmented parent-level entropy account.  Once the next parent
mean is known, the entropy of that parent mean is known, but the decomposition
into subcell Jensen gap, entropy-flux divergence, and residual dissipation is
additional information.

Let \(\eta\) be a convex entropy on a convex state domain.  Define the parent
mean and subcell Jensen gap in parent \(K\) at time level \(m\) by
\[
  \bar U_K^m=(RU^m)_K,
  \qquad
  \delta_K^m
  =
  \sum_{i\subset K}\frac{h_i}{H_K}\eta(U_i^m)
  -
  \eta(\bar U_K^m).
\]
By Jensen's inequality, \(\delta_K^m\ge0\).

\begin{theorem}[Augmented parent entropy balance]
\label{thm:entropybalance}
Suppose the fine cells satisfy the exact entropy-amount balance
\[
  h_i\eta(U_i^{n+1})+\mathcal D_i
  =
  h_i\eta(U_i^n)-(\Psi_{i+1/2}-\Psi_{i-1/2}),
  \qquad
  \mathcal D_i\ge0,
\]
for a declared numerical entropy flux \(\Psi_{i+1/2}\).  Let
\(\Psi_{K,+}\) and \(\Psi_{K,-}\) denote the realized integrated entropy
fluxes through the right and left parent interfaces.  Then
\[
  \eta(\bar U_K^{n+1})
  +
  \delta_K^{n+1}
  +
  \frac{\mathcal D_K}{H_K}
  =
  \eta(\bar U_K^{n})
  +
  \delta_K^{n}
  -
  \frac{\Psi_{K,+}-\Psi_{K,-}}{H_K},
\]
where
\[
  \mathcal D_K=\sum_{i\subset K}\mathcal D_i.
\]
\end{theorem}

\begin{proof}
Sum the fine entropy balances over all children in parent \(K\).  Internal
entropy fluxes telescope.  The identity
\[
  \sum_{i\subset K}h_i\eta(U_i^m)
  =
  H_K\bigl[\eta(\bar U_K^m)+\delta_K^m\bigr]
\]
then gives the displayed parent balance after division by \(H_K\).
\end{proof}

This theorem is conditional.  It assumes a cellwise entropy balance with a
declared numerical entropy flux and nonnegative residuals.  It does not claim
that every finite-volume or Runge--Kutta method automatically has such a
balance.  Different entropy-flux conventions can redistribute entropy flux
and residual dissipation, so the information separation below is for the
declared entropy pair and declared numerical entropy flux.

\subsection{A finite entropy-blind Burgers--Rusanov family}

Consider inviscid Burgers' equation
\[
  u_t+\left(\frac{u^2}{2}\right)_x=0
\]
with quadratic entropy and entropy flux
\[
  \eta(u)=\frac12u^2,
  \qquad
  q(u)=\frac13u^3.
\]
Use forward Euler with \(\Delta t=1/2\), unit fine cells, and the Rusanov
conservative numerical flux with \(\alpha=1\),
\[
  \widehat f(u_L,u_R)
  =
  \frac12\bigl(f(u_L)+f(u_R)\bigr)
  -
  \frac12(u_R-u_L).
\]
We pair it with the declared Rusanov-type numerical entropy flux
\[
  \widehat q(u_L,u_R)
  =
  \frac12\bigl(q(u_L)+q(u_R)\bigr)
  -
  \frac12\bigl(\eta(u_R)-\eta(u_L)\bigr).
\]
The corresponding integrated conservative and entropy fluxes are
\[
  \Phi_{i+1/2}=\Delta t\,\widehat f(u_i,u_{i+1}),
  \qquad
  \Psi_{i+1/2}=\Delta t\,\widehat q(u_i,u_{i+1}).
\]

\begin{lemma}[Local entropy residuals]
\label{lem:local-entropy-residuals}
Let \(a,b\in\{-1,+1\}\), and consider the four-cell local state
\[
  (a,-a,-b,b),
\]
with compatible exterior states \(a\) on the left and \(b\) on the right.
After one Rusanov forward-Euler step with \(\Delta t=1/2\), the four updated
values are
\[
  \left(
  \frac a2,\,
  -\frac{a+b}{4},\,
  -\frac{a+b}{4},\,
  \frac b2
  \right).
\]
The residuals
\[
  \mathcal D_i
  =
  \eta(u_i^n)-\eta(u_i^{n+1})
  -
  \bigl(\Psi_{i+1/2}-\Psi_{i-1/2}\bigr)
\]
are nonnegative and are listed in \Cref{tab:local-entropy-residuals}.
\end{lemma}

\begin{table}[htbp]
\centering
\caption{Local entropy residuals for the four-cell Burgers--Rusanov stencil
\((a,-a,-b,b)\), with \(a,b\in\{-1,+1\}\).  All entries are nonnegative,
which verifies the cellwise entropy-balance hypothesis used in
\Cref{thm:entropybalance} for every local branch of the sign family.}
\label{tab:local-entropy-residuals}
\begin{tabular}{@{}ccccc@{}}
\toprule
\((a,b)\) & \(\mathcal D_0\) & \(\mathcal D_1\) & \(\mathcal D_2\) & \(\mathcal D_3\) \\
\midrule
\((1,1)\)     & \(13/24\) & \(13/24\) & \(5/24\)  & \(5/24\)  \\
\((1,-1)\)    & \(13/24\) & \(1/2\)   & \(1/2\)   & \(13/24\) \\
\((-1,1)\)    & \(5/24\)  & \(1/2\)   & \(1/2\)   & \(5/24\)  \\
\((-1,-1)\)   & \(5/24\)  & \(5/24\)  & \(13/24\) & \(13/24\) \\
\bottomrule
\end{tabular}
\end{table}

\begin{proof}
For \(u_L,u_R\in\{-1,+1\}\),
\[
  \widehat f(u_L,u_R)
  =
  \frac12+\frac12(u_L-u_R).
\]
Substitution into the local forward-Euler update gives the displayed updated
values.  The numerical entropy flux satisfies
\[
  \widehat q(u_L,u_R)
  =
  \frac{u_L^3+u_R^3}{6}
  -
  \frac{u_R^2-u_L^2}{4}.
\]
For sign states, \(u_L^2=u_R^2=1\), so the second term vanishes and
\[
  \Psi_{i+1/2}
  =
  \frac{u_i+u_{i+1}}{12}.
\]
Substitution into the definition of \(\mathcal D_i\) gives the four rows in
\Cref{tab:local-entropy-residuals}.
\end{proof}

\begin{theorem}[Conservative blindness to entropy-accounting data]
\label{thm:entropyblind}
For every \(P\ge2\), let
\[
  \sigma=(\sigma_0,\ldots,\sigma_{P-1})\in\{-1,+1\}^P.
\]
In parent \(K\), set the four child values to
\[
  (\sigma_{K-1},-\sigma_{K-1},-\sigma_K,\sigma_K),
\]
with periodic indexing.  Then all \(2^P\) branches have identical parent
means, identical conservative parent-interface flux registers, and identical
next parent means.  However, their raw entropy-interface registers are
\[
  \Psi_K(\sigma)=\frac{\sigma_K}{6},
\]
which distinguish all branches.  The entropy-flux part of the periodic
parent entropy balance depends on the cyclic divergence
\[
  \beta_K(\sigma)
  =
  \Psi_K(\sigma)-\Psi_{K-1}(\sigma)
  =
  \frac{\sigma_K-\sigma_{K-1}}{6},
\]
which takes \(2^P-1\) distinct values.  Under the uniform branch law,
\[
  \mathbb P(\beta=0)=2^{1-P}.
\]
Moreover, for deterministic prediction from the conservative channel under
squared Euclidean loss,
\[
  \inf_{\widehat\Psi}\mathbb E\|\Psi-\widehat\Psi\|_2^2
  =
  \frac{P}{36},
  \qquad
  \inf_{\widehat\beta}\mathbb E\|\beta-\widehat\beta\|_2^2
  =
  \frac{P}{18}.
\]
\end{theorem}

\begin{proof}
Each parent has zero mean.  The right interface of parent \(K\) has equal
left and right states \(\sigma_K\).  Hence the conservative integrated flux
through every parent interface is
\[
  \Delta t\, f(\sigma_K)=\frac12\cdot\frac12=\frac14,
\]
independent of \(\sigma\).  All conservative flux differences vanish, so the
next parent means are identical.

At the same equal interface states,
\[
  \Psi_K=\Delta t\,q(\sigma_K)=\frac12\cdot\frac{\sigma_K}{3}
  =
  \frac{\sigma_K}{6}.
\]
Thus \(\Psi(\sigma)=\sigma/6\) distinguishes all \(2^P\) sign branches.

Now consider \(\beta_K=(\sigma_K-\sigma_{K-1})/6\).  If
\(\beta=0\), then every neighboring pair is equal, so \(\sigma\) is one of
the two constant sign vectors.  If \(\beta\ne0\), the cyclic differences
determine the entire sign vector: at each edge, a zero difference means the
sign is unchanged, while a nonzero difference gives the direction of the
sign flip.  Hence all nonzero divergence patterns are distinct, and the
number of divergence values is \(2^P-1\).  Under the uniform branch law,
\[
  \mathbb P(\beta=0)=\frac{2}{2^P}=2^{1-P}.
\]

For squared-error prediction from the conservative channel, the input channel
is constant over all branches.  Under the uniform sign law,
\[
  \mathbb E[\Psi]=0,
  \qquad
  \mathbb E[\beta]=0.
\]
Thus the optimal deterministic predictors are zero.  Since
\[
  \|\Psi\|_2^2=\sum_{K=0}^{P-1}\frac{\sigma_K^2}{36}=\frac{P}{36},
\]
the first error floor follows.  For the divergence,
\[
  \mathbb E\|\beta\|_2^2
  =
  \frac1{36}
  \sum_{K=0}^{P-1}\mathbb E[(\sigma_K-\sigma_{K-1})^2].
\]
Independent neighboring signs differ with probability \(1/2\), so
\[
  \mathbb E[(\sigma_K-\sigma_{K-1})^2]=2,
\]
and the second error floor is \(P/18\).  Finally,
\Cref{lem:local-entropy-residuals} verifies the nonnegative cellwise
residuals required by \Cref{thm:entropybalance}.
\end{proof}

\begin{table}[htbp]
\centering
\small
\caption{Entropy-accounting sign-family counts and error floors.  For each
\(P\), all \(2^P\) branches share the conservative channel.  The raw entropy
register distinguishes all branches, while the gauge-invariant cyclic
divergence has \(2^P-1\) values because the two constant sign states share
zero divergence.  The final columns give the deterministic squared-error
floors under the uniform branch law.}
\label{tab:entropy-counts}
\begin{tabular}{@{}rrrrrrr@{}}
\toprule
\(P\) & branches & \(\Psi\) values & \(\beta\) values
& \(\mathbb P(\beta=0)\)
& \(\mathbb E\|\Psi\|_2^2\)
& \(\mathbb E\|\beta\|_2^2\) \\
\midrule
\(2\)  & \(4\)      & \(4\)      & \(3\)      & \(1/2\)       & \(1/18\) & \(1/9\) \\
\(4\)  & \(16\)     & \(16\)     & \(15\)     & \(1/8\)       & \(1/9\)  & \(2/9\) \\
\(8\)  & \(256\)    & \(256\)    & \(255\)    & \(1/128\)     & \(2/9\)  & \(4/9\) \\
\(16\) & \(65{,}536\) & \(65{,}536\) & \(65{,}535\) & \(1/32{,}768\) & \(4/9\)  & \(8/9\) \\
\bottomrule
\end{tabular}
\end{table}

\Cref{tab:local-entropy-residuals} connects the sign family to the cellwise
entropy-balance hypothesis.  \Cref{tab:entropy-counts} summarizes the global
information separation.  The conservative channel is identical across all
branches, but the entropy register and its gauge-invariant divergence are
not.  The error floors are not distribution-free information-theoretic
bounds; they are optimal deterministic squared-error risks for the uniform
branch law and the declared entropy targets. The table entries and finite certificate checks are also generated by the
accompanying Python scripts using exact rational arithmetic via
\texttt{fractions.Fraction}.  Thus the reported finite identities, residual
values, branch counts, and error floors do not depend on floating-point
roundoff.  The scripts serve as reproducibility certificates for the finite
examples; the mathematical claims are proved analytically above.

\section{Discussion}
\label{sec:discussion}

The results separate one-step information sufficiency from recursive
predictive memory.  Recursive memory, conditioning, and dissipative decay are
studied in the companion paper \cite{PolemitisChristakisDrikakis2026FVerasure}.
Here the focus is narrower: which information channels determine one
completed transition and its entropy account?

The SSPRK(2,2) example shows that initial traces are integrator-dependent
information objects.  Forward Euler uses only the initial trace pair at a
parent interface.  A multistage method can generate later stage traces from
hidden interior values, so the correct ex post object is the realized
stage-integrated flux register.

The entropy construction shows that conservative closure and entropy
accounting are separate contracts.  The conservative channel in
\Cref{thm:entropyblind} determines next parent means but not the augmented
entropy account.  In learned, reduced, and multiscale finite-volume models
\cite{BarSinai:2019datadriven,KochkovEtAl2021,AgdesteinSanderse2025,
SanderseEtAl2025}, this suggests that conservation, entropy accounting,
stage information, and prediction horizon should be specified separately.

\subsection{Limitations and future work}
\label{sec:limitations}

The analysis is intentionally narrow.  The SSPRK(2,2) obstruction is a scalar
periodic upwind construction, and the entropy-accounting separation uses a
particular Burgers--Rusanov entropy-flux convention.  The results do not
claim a universal lower bound for compressible Euler systems, high-order
weighted essentially non-oscillatory reconstructions, nonperiodic grids, or
learned closures.

Natural extensions include a general explicit Runge--Kutta theorem
characterizing how far interior information can enter stage traces, and a
more general entropy-separation theorem identifying flux and entropy pairs
that admit conservative collisions with entropy-accounting separation.
Another direction is to test the same channel distinctions at adaptive
mesh-refinement interfaces and for system examples.

\section{Conclusions}
\label{sec:conclusions}

A one-step finite-volume closure exists only for declared information that
makes the restricted fine update single valued.  Complete ordered initial
traces close a forward-Euler finite-volume transition with a deterministic
two-state numerical flux, but they do not generally close an SSPRK(2,2)/Heun
transition.  The realized stage-integrated interface flux remains sufficient
for the completed transition, but it is an ex post register.

Conservative one-step closure is also not entropy-accounting closure.  For a
declared Burgers--Rusanov entropy account, the sign family constructed here
has identical conservative parent means, conservative interface registers,
and next parent means across all branches, while the entropy register and its
cyclic divergence vary.  Under the uniform branch law, the conservative
channel leaves nonzero deterministic squared-error floors for the entropy
targets.

The practical message is that coarse, reduced, and learned finite-volume
models should specify their retained averages, interface traces, stage
registers, entropy-side registers, and prediction horizon before claiming
pathwise closure.

\section*{Funding}

This work received no external funding.

\section*{Declaration of competing interest}

The authors declare that they have no known competing financial interests or
personal relationships that could have appeared to influence the work
reported in this paper.

\section*{CRediT authorship contribution statement}

Antonis Polemitis: Conceptualization, Methodology, Verification,
Writing--original draft, Writing--review and editing.  Nicholas Christakis:
Conceptualization, Methodology, Formal analysis, Software, Verification,
Writing--original draft, Writing--review and editing.  Dimitris Drikakis:
Conceptualization, Methodology, Verification, Writing--review and editing.

\section*{Data and code availability}

The complete reproducibility package for this study is publicly
available in the Certified Simulation repository of the University of
Nicosia, \url{https://github.com/UniversityOfNicosia/certified-simulation},
under \texttt{papers/initial-trace-entropy-separation}.  The package
contains the scripts that generate the computational certificates, CSV
files, and LaTeX table fragments reported in the manuscript, together
with an independent certificate suite that replays the paper's exact
results in exact rational arithmetic. The table entries and finite certificate checks are also generated by the accompanying Python scripts using exact rational arithmetic via \texttt{fractions.Fraction}.  Thus the reported finite identities, residual values, branch counts, and error floors do not depend on floating-point roundoff.  The scripts serve as reproducibility certificates for the finite examples; the mathematical claims are proved analytically above.  No external CFD solver or proprietary software is required.




\FloatBarrier
\bibliographystyle{elsarticle-num} 

\bibliography{bibliography}


\end{document}